\documentclass{article}
\usepackage[dvips]{graphicx}
\usepackage{graphicx}
\usepackage{amsmath, mathrsfs, wasysym}% txfonts}
\usepackage{amssymb}
\usepackage{ifsym}
\usepackage{textcomp}
\usepackage{stmaryrd}
\usepackage{marvosym}
\usepackage{latexsym, amsmath, amssymb}
\usepackage{color}
\usepackage{multirow}

\newtheorem{theorem}{Theorem}[section]

\newtheorem{definition}[theorem]{Definition}

\newtheorem{prop}[theorem]{Proposition}
\newtheorem{remark}[theorem]{Remark}

\numberwithin{equation}{section}

\newcommand{\abs}[1]{\lvert#1\rvert}

\newcommand{\x} {{\bf{x}}}
\newcommand{\OL} {\mathcal{L}}
\newcommand{\ul}[1]{\underline{\mathbf{{#1}}}}
\newcommand{\ull}[1]{\underline{\underline{\mathbf{{#1}}}}}
\newcommand{\A}{\mathcal{A}}
\newcommand{\K}{\kappa}\newcommand{\p}{\partial}
\newcommand{\na}{{\nabla_{x}}}

\newcommand{\sg}{\underline{\sigma}}
\newcommand{\st}{\underline\varepsilon}
\newcommand{\bx}{{\mathbf{x}}}
\newcommand{\Gc}{\underline{\mathbf{\Gamma}}^{c}}
\newcommand{\Gv}{\underline{\mathbf{\Gamma}}^{v}}
\newcommand{\G}{\underline{\mathbf{\Gamma}}}

\newcommand{\eqnref}[1]{(\ref {#1})}
\newcommand{\href}[1]{(H.\ref {#1})}
\newcommand{\RR}{\mathbb{R}}
\newcommand{\be}{\begin{equation}}
\newcommand{\ee}{\end{equation}}
\newcommand{\dis}{\displaystyle}
\newcommand{\E}{{\bf{\underline{E}}}}
\newcommand{\D}{\mathbf{D}}
\newcommand{\bPsi}{{\overline{\Psi}}}
\newcommand{\hr}{{\widehat{\bf{r}}}}
\newcommand{\hR}{{\widehat{\bf{R}}}}
\newcommand{\bu}{{\bf{u}}}
\newcommand{\tphi}{{\widetilde{\phi}}}
\newcommand{\tPhi}{{\widetilde{\Phi}}}
\newcommand{\bPhi}{{\overline{\Phi}}}
\newcommand{\bef}{{{\bf{e}}}}
\newcommand{\bn}{{\bf{n}}}
\newcommand{\bG}{\underline{{\bf{G}}}}
\newcommand{\bF}{{\mathbf{F}}}

\title{Some Anisotropic Viscoelastic Green Functions\footnote{Submitted  to Contemporary Mathematics, American Mathematical Society, Providence.}}
\author{Elie Bretin\thanks{Centre de
Math{\'e}matiques Appliqu{\'e}es, CNRS UMR 7641, Ecole
Polytechnique, 91128 Palaiseau, France
(bretin@cmap.polytechnique.fr,
wahab@cmap.polytechnique.fr).} \and Abdul Wahab\footnotemark[2]}

\begin{document}
%%%
%%% CMAP %%%
%\pagedegarde{Some Anisotropic Viscoelastic Green Functions} {Elie Bretin, Abdul Wahab} {708} {February  2011}
%%%
%%%

\maketitle

\begin{abstract}
In this paper, we compute the closed form expressions of elastodynamic 
Green functions for three different viscoelastic media with
simple type of anisotropy. 
We follow Burridge \emph{et al.}
[Proc. Royal Soc. of London. 440(1910): (1993)] to express 
unknown Green function in terms of three scalar functions $\phi_i$, by 
using the spectral decomposition of the Christoffel tensor
associated with the medium. 
The problem of computing Green function is, thus reduced to
the resolution of three scalar wave equations satisfied by $\phi_i$, 
and subsequent equations with $\phi_i$ as source terms.
To describe viscosity effects, we choose an empirical power law model  
which becomes well known Voigt model for quadratic frequency losses.
\end{abstract}

\maketitle

\section{Introduction}

\par
Numerous applications in biomedical 
imaging \cite{AmmariBioMed, AmmariRSI}, seismology 
\cite{AkiRich, Carcoine2007Wave}, exploration geophysics \cite{Helbig, HelbigThomsen}, 
 material sciences \cite{ShapeOpt, AmmariPol} and engineering sciences 
\cite{Achenbach, BinSingh, lekhnitskii1964theory} have fueled research 
and development in theory of elasticity. 
Elastic properties and 
attributes have gained interest in the recent decades as a diagnostic tool for non-invasive imaging \cite{Greenleaf, Sarvazyan}. Their high correlation
with the pathology and the underlying structure of soft tissues 
has inspired many investigations in biomedical imaging
and led to many interesting mathematical problems
\cite{ AmmariDirEast, AmmariMRElasto, AmmariOpti, 
 AmmariEffVis, AmmariSepration, Bercoff, Sinkus00, Sinkus05}.

\par
Biological materials are often assumed to be isotropic and inviscid 
with respect to elastic deformation. However, several recent studies indicate
that many soft tissues exhibit anisotropic and viscoelastic behavior 
\cite{GenCathChafFink, OidaKang, Sinkus00, Sinkus05, NamaniBaylay09, YoonKatz}. 
Sinkus \emph{et al.}  have inferred in \cite{Sinkus00} that breast tumor 
tends to be anisotropic, while Weaver \emph{et al.} \cite{Weaver} have provided
 an evidence that even non cancerous breast tissue is anisotropic.
White matter in brain \cite{NamaniBaylay09} and cortical bones \cite{YoonKatz}
 also exhibit similar behavior. Moreover, it has been observed that the shear
velocities parallel  and orthogonal to the fiber direction in forearm 
\cite{OidaKang} and biceps \cite{GenCathChafFink} are different. This
 indicates that the skeletal muscles with directional structure
are actually anisotropic.
Thus, an assumption of isotropy can lead to erroneous 
forward-modelled wave synthetics, while an estimation of viscosity effects
 can be very useful in 
characterization and identification of anomaly \cite{Bercoff}.

\par
A possible approach to handle viscosity effects on image reconstruction
 has been proposed in \cite{ELA11} using stationary phase theorem. It is shown
that the ideal Green function
 (in an inviscid regime) can be approximated from the viscous one
by solving an ordinary differential equation. Once the ideal Green function 
is known one can identify a possible anomaly using imaging algorithms such as
time reversal, back-propagation, Kirchhoff migration or MUSIC 
\cite{AmmariDirEast, AmmariTransElas, AmmariRSI, AmmariBioMed}.  
One can also find the elastic moduli of the anomaly using the 
asymptotic formalism and reconstructing a certain polarization tensor in the 
far field
\cite{AmmariMRElasto, AmmariTransElas, AmmariPol, AmmariExp}.

\par
The importance of Green function stems from its role as a
tool for the numerical and asymptotic techniques in biomedical imaging. 
Many inverse problems
involving the estimation and acquisition of elastic parameters become
tractable once the associated Green function is computed 
\cite{AmmariBound, AmmariDirEast, AmmariTransElas, ELA11}.
Several attempts have been made to compute 
Green functions in purely elastic  and/or isotropic regime. 
(See \emph{e.g.} \cite{ELA11, Bercoff, Burridge93, Carcoine2007Wave, Payton,
 VavAsymtGreen07, vavrycuk2001exact, vavrycukWeak} 
and references therein). However, it is not possible to give a closed 
form expression for general anisotropic Green functions without certain restrictions on the media. 
In this work, we provide anisotropic viscoelastic Green function in closed form
for three particular anisotropic media.

\par
The elastodynamic Green function in isotropic media is calculated 
by separating wave modes using Helmholtz 
decomposition of the elastic wavefield \cite{AkiRich, ELA11, Bercoff}.
Unfortunately, this simple approach does not work in anisotropic 
media, where three different waves propagate with different phase velocities
 and polarization directions 
\cite{Carcoine2007Wave, BinSingh, cerveny2008RayTheory}. A polarization direction of 
quasi-longitudinal wave that differs from that of wave vector, impedes
Helmholtz decomposition to completely separate wave modes \cite{Dellinger}.

\par
The phase velocities and polarization vectors are the eigenvalues
and  eigenvectors of the Christoffel tensor $\G$ associated with the 
medium. 
So, the wavefield can always be decomposed using the spectral basis of 
 $\G$. Based on this observation, Burridge \emph{et al.} \cite{Burridge93} 
proposed a new approach to calculate elastodynamic Green functions.
Their approach consists of finding the eigenvalues and 
eigenvectors of Christoffel tensor $\G(\na)$ using the duality between
algebraic and differential objects.
Therefore it is possible to express the Green function $\bG$ 
in terms of three scalar functions $\phi_i$ satisfying partial differential
 equations with constant coefficients.
Thus the problem of computing $\bG$ reduces to the 
 resolution of three differential equations for $\phi_i$
and of three subsequent equations (which may or may not 
be differential equations) with $\phi_i$  as source terms. See 
\cite{Burridge93} for more details.

\par
Finding the closed form expressions of the eigenvalues of 
Christoffel tensor $\G$ is usually not 
so trivial because its characteristic equation is 
a polynomial of degree six in the components of its argument vector. 
However, with some restrictions on the material, roots of the
characteristic equation can be given \cite{Payton}. In this article, 
we consider three different media for which not only the explicite expressions of the 
eigenvalues of $\G$ are known \cite{Burridge93, vavrycuk2001exact}, but they are also 
quadratic homogeneous forms, in the components of the 
argument vector. As a consequence, 
equations satisfied by $\phi_i$ become scalar wave equations.
Following Burridge \emph{et al.}  \cite{Burridge93}, 
we find the viscoelastic  Green functions
 for each medium. It is important to note that the
 elastodynamic Green function in a purely elastic regime, for the media under consideration, are well known \cite{vavrycuk2001exact, Burridge93}.
Also, the expression of the Green function for viscoelastic isotropic medium,  which is computed as a special case, matches the one provided 
in \cite{ELA11}.

\par
It has been shown in \cite{Catheline} that Voigt model is well adopted to 
describe the viscosity response of many soft tissues to low frequency 
excitations. In this work, we consider a more general model proposed by Szabo and Wu in \cite{SzaboWu00}, which 
describes an empirical power law  behavior of many viscoelastic
materials including human myocardium. This model is based on a time-domain
 statement of causality \cite{SzaboCausal94, titchmarsh} and reduces
to Voigt model for the specific case of quadratic frequency losses.

\par
We provide some mathematical notions, theme and the outlines of 
the article in the next section.

\section{Mathematical Context and Paper Outlines}

\subsection{Viscoelastic Wave Equation}

\par
Consider an open subset $\Omega$ of $\RR^3$, filled 
with a homogeneous anisotropic viscoelastic material. Let 
$$
\bu (\x,t):\Omega\times\RR^+\to\RR^3
$$ 
be the displacement field at time $t$ of the material 
particle at position $\x\in\Omega$ and $\na\bu(\x,t)$ be its 
gradient. 

Under the assumptions of linearity and small perturbations, 
we define the order two strain tensor 
 by
\be 
\st:(\x,t)\in\Omega\times\RR^+\longmapsto\frac{1}{2}
\left(\na\bu+\na\bu ^T\right)(\x,t),
\ee 
where the superscript $T$ indicates a transpose operation.

\par
Let $\ull C\in\OL_s^2(\RR^3)$ and $\ull V\in\OL_s^2(\RR^3)$ be
the stiffness and viscosity tensors of the material respectively. Here  $\OL^2_s(\RR^3)$ is the space of symmetric tensors of order four.
These tensors are assumed to be positive definite, i.e.  there exists a constant $\delta>0$ such that
$$
(\ull C : \ul \xi):\ul \xi\geq\delta\abs{\ul\xi}^2  \quad \text{and} \quad  (\ull V : \ul \xi):\ul \xi\geq\delta\abs{\ul\xi}^2, \quad\forall \ul\xi\in\OL_s(\RR^d),
$$
where $\OL_s(\RR^3)$ denotes the space of symmetric tensors of order two.

\par
The generalized Hooke's Law \cite{SzaboWu00} for power law media 
 states that the stress distribution 
$$\sg:\Omega\times\RR^+\to\OL_s(\RR^3)$$
produced by deformation $\st$, satisfies:
\be 
\sg=\ull C:\st+\ull V:\A[\st]\label{HookeLaw}
\ee
where $\A$ is a causal operator defined as 
\begin{equation} \label{defA}
\A[\varphi] = \left|
\begin{array}{ll}
\dis-(-1)^{\gamma/2}\frac{\partial^{\gamma-1}\varphi}{\partial t^{\gamma-1}}
&\gamma\text{  is an even integer},
\\
\\
\dis\frac{2}{\pi}(\gamma-1)!(-1)^{(\gamma+1)/2}
\frac{H(t)}{t^\gamma}*_t \varphi
&\gamma\text{ is an odd integer},
\\
\\
\dis-\frac{2}{\pi}{\Gamma}(\gamma)\sin(\gamma\pi/2)
\frac{H(t)}{|t|^\gamma}*_t \varphi
&\gamma\text{ is a non integer}.
\end{array}
\right.
\end{equation}
Note that by convention,  
$$
\A[\bu]_i=\A[u_i]
\quad\text{and}\quad 
\A[\st]_{ij}=\A[\varepsilon_{ij}]\quad 1\leq i,j\leq 3. 
$$
Here $H(t)$ is the Heaviside function, $\Gamma$ is the gamma function
 and $*_t$ represents convolution with respect to variable $t$. 
See \cite{ AlekseevRybak, Caputo, SzaboWu00, SzaboCausal94, 
titchmarsh} for 
comprehensive details and discussion on fractional attenuation models,
causality and the loss operator $\A$.

\par
The viscoelastic  wave equation satisfied by the displacement field $\bu(\x,t)$ reads now 
\begin{eqnarray*}
 \rho\frac{\partial^2 \bu}{\partial t^2}-\mathbf{F} &=& \na\cdot\sg = \na\cdot \left(\ull C:\st+\ull V:\A[\st] \right),
\end{eqnarray*}
where  $\bF(\x,t)$ is the applied force and $\rho$ is the  
density (supposed to be constant) of the material.

\begin{remark}
For quadratic frequency losses, i.e, when $\gamma = 2$,  operator $\A$ reduces
 to a first order time derivative. Therefore, power-law attenuation model
turns out to be the Voigt model in this case.
\end{remark}

%%%%%%%%%%%%%%%%
\subsection{Spectral decomposition by Christoffel tensors} 
%%%%%%%%%%%%%%%%%%

\par
We introduce now the Christoffel tensors  $\Gc, \Gv:\RR^3\to\OL_s(\RR^3)$   associated respectively with $\ull C$ and $\ull V$ defined by: 
$$
 \Gamma^c_{ij}(\bn)=\sum_{k,l=1}^3C_{kilj}n_kn_j, \quad  \Gamma^v_{ij}(\bn)=\sum_{k,l=1}^3V_{kilj}n_kn_j, \quad \quad  \quad\forall\bn\in\RR^3,\quad 1\leq i,j\leq 3.
$$
Remark that the viscoelastic wave equation can be rewritten in terms of  Christoffel tensors as : 
\begin{eqnarray}
\rho\frac{\partial^2 \bu}{\partial t^2}-\mathbf{F}  
= \Gc[\na] \bu + \Gv[\na] \A[\bu].\label{eq:ElastWave}
\end{eqnarray}
Note that $\Gc$ and $\Gv$ are symmetric and positive definite as  $\ull C$ and $\ull V$  are already symmetric positive definite.

\par
Let $L_i^c$ be the eigenvalues and $\D_i^c$ be the associated
 eigenvectors of $\Gc$ for $i=1,2,3$.
We define the quantities $M_i^c$ and $\E_i^c$ by 
\be 
M_i^c=\D_i^c\cdot \D_i^c ,
\quad\text{and}\quad
\E_i^c= (M_i^c)^{-1}\D_i^c\varotimes \D_i^c.
\ee

\par
As $\Gc$ is symmetric, the eigenvectors $\D_i^c$ are orthogonal and the spectral decomposition of the Christoffel  tensor $\Gc$ can be given as:
\be 
\Gc=\sum_{i=1}^3L^c_i\E^c_i\quad\text{with}\quad \ul I=
\sum_{i=1}^3\E^c_i\label{SpectGamma} 
\ee 
where $\ul I \in \OL_s(\RR^3)$ is the identity tensor.

\par
Similarly,  consider  $\Gv$ the Christoffel tensor associated  with  $\ull V$ and define the quantities 
 $L_i^v$,  $\D_i^v$, $M_i^v$ and $\E_i^v$ such as 
\be 
\Gv=\sum_{i=1}^3L^v_i\E^v_i\quad\text{with}\quad \ul I=
\sum_{i=1}^3\E^v_i\label{SpectGammav}.
\ee 

\par
We assume that the tensors $\Gc$ and $\Gv$ have the same structure in the sense that the eigenvectors  $\D_i^c$ and $\D_i^v$ are  equal. (See Remark \ref{RemD}). In the sequel we use $\D$ instead of $\D^c$ or  $\D^v$ and similar for $\E$ and $M$, by abuse of notation.

\subsection{Paper Outline}

\par
The aim of this work is to compute the elastodynamic Green function $\bG$ associated 
to viscoelastic wave equation \eqnref{eq:ElastWave}. More precisely,  $\bG$ is the solution of the equation 
\be 
\left( \Gc[\na]\bG (\bx,t) + \Gv[\na]\A[\bG](\bx,t) \right) 
-\rho\dis\frac{\p^2\bG(\bx,t)}{\p t^2}
=\delta(t)\delta(\bx){\ul I},  
\label{eq:Green}
\ee 
The idea is to use the spectral decomposition of $\bG$ of the form 
\be 
\bG=\dis\sum_{i=1}^3\E_i(\na)\phi_i
=\dis\sum_{i=1}^3 (\D_i\otimes\D_i)M_i^{-1}\phi_i,\label{PhiSpectSum}
\ee
where $\phi_i$ are three scalar functions satisfying  
\be 
%\begin{array}{ll}
\left( L^c_i(\na)\phi_i + L^v_i(\na) \A[\phi_i] \right) 
-\rho\dis\frac{\p^2\phi_i}{\p t^2}=\delta(t)\delta(\bx)  
%\end{array}
%\right.
\label{eq:PHI}
\ee
(See Appendix  \ref{appendix_phi} for more details 
about this decomposition.)

Therefore, to obtain an expression of  $\bG$, we need to: 
\begin{itemize}
 \item[1-] solve  three partial differential equations \eqref{eq:PHI} in $\phi_i$
\item[2-] subsequent  equations 
\be
\psi_i=M_i^{-1} \phi_i\label{PotentialProblem}
\ee 
\item[3-] and calculate second order derivatives of  $\psi_i$ to compute  
$$  (\D_i\otimes\D_i) \psi_i$$
\end{itemize}

\par 

In the following Section, we give simple examples of anisotropic media 
which satisfy some restrictive properties and  assumptions 
(see Subsection \ref{Assumptions}) defining the limits of our approach. 
In Section \ref{WaveProblem}, we derive the solutions 
 $\phi_i$ of equations \eqref{eq:PHI}.
In Section  \ref{section_potential}, we give an explicite 
resolution of $\psi_i=M_i^{-1} \phi_i$ and 
$  (\D_i\otimes\D_i) \psi_i$.
Finally, in the last section, we compute the Green function
for three simple anisotropic media.

\section{Some Simple Anisotropic Viscoelastic Media}\label{sec:Media}

\par 
In this section, we present three viscoelastic media with 
simple type of anisotropy. We also describe some important properties of
the media and our basic assumptions in this article.

\begin{definition}
We will call a tensor $\ul c=\left(c_{mn}\right)\in\OL_s(\RR^6)$ the Voigt 
representation of an order four tensor $\ull C\in \OL^2_s(\RR^3)$ if 
$$
c_{mn}=c_{p(i,j)p(k,l)}=C_{ijkl}\quad 1\leq i,j,k,l\leq 3 
$$
where
\begin{equation*}
p(i,i)=i,\quad p(i,j)=p(j,i), \quad p(2,3)=4,\quad p(1,3)=5,\quad p(1,2)=6.
\end{equation*}
We will use $\ul c$ and $\ul v$ for the Voigt representations of
stiffness tensor $\ull C$ and viscosity  tensor $\ull V$ respectively.
\end{definition}

\par
We will let tensors $\ul c$ and $\ul v$ to have a same structure.
For each media, the expressions for 
$\Gc$, $L_i^c(\na)$,  $\D_i^c(\na)$ and $M_i^c(\na)$ are provided
\cite{Burridge93, vavrycuk2001exact}.
Throughout this section, $\mu_{pq}$ will assume the value $c_{pq}$ 
for $\ul c$ and $v_{pq}$ for $\ul v$ where the subscripts 
$p,q\in\{1,2,\cdots,6\}$. Moreover, we assume that 
the axes of material are identical with the Cartesian coordinate
axes $\bef_1, \bef_2$ and $\bef_3$ and $\p_i=\dis\frac{\p}{\p x_i}$.

\subsection{Medium I}

\par 
The first medium for which we present a closed form elastodynamic Green
function is an orthorhombic 
medium with the tensors $\ul c$ and $\ul v$ of the form:
$$
\begin{pmatrix}
\mu_{11} & -\mu_{66} & -\mu_{55} & 0      & 0      & 0\\
-\mu_{66} & \mu_{22} & -\mu_{44} & 0      & 0      & 0\\
-\mu_{55} & -\mu_{44} & \mu_{33} & 0      & 0      & 0 \\
0      & 0      & 0      & \mu_{44} & 0      & 0\\ 
0      & 0      & 0      & 0      & \mu_{55} & 0\\
0      & 0      & 0      & 0      & 0      & \mu_{66}
\end{pmatrix}
$$
The Christoffel tensor is given by
$$
\Gc=
\begin{pmatrix}
c_{11}\p_1^2+c_{66}\p_2^2+c_{55}\p_3^2 & 0 & 0\\
0 & c_{66}\p_1^2+c_{22}\p_2^2+c_{44}\p_3^2 & 0\\
0 & 0 & c_{55}\p_1^2+c_{44}\p_2^2+c_{33}\p_3^2\\
\end{pmatrix}
$$
Its  eigenvalues $L_i^c(\na)$ and the associated eigenvectors $\D_i^c(\na)$  are:
\begin{eqnarray*}
L_1^c(\na) = c_{11}\p_1^2+c_{66}\p_2^2+c_{55}\p_3^2\\
L_2^c(\na) = c_{66}\p_1^2+c_{22}\p_2^2+c_{44}\p_3^2\\
L_3^c(\na) = c_{55}\p_1^2+c_{44}\p_2^2+c_{33}\p_3^2\\
\D_i^c = {\bef}_i\quad\text{with}\quad M_i^c=1\quad\forall i=1,2,3\\
\end{eqnarray*}

\subsection{Medium II}

\par
The second medium which we consider is a transversely 
isotropic medium having symmetry axis along $\bef_3$ and 
defined by the stiffness and the viscosity tensors $\ul c$ and $\ul v$ 
of the form:
$$
\begin{pmatrix}
\mu_{11} & \mu_{12} & -\mu_{44} & 0      & 0      & 0\\
\mu_{12} & \mu_{11} & -\mu_{44} & 0      & 0      & 0\\
-\mu_{44} & -\mu_{44} & \mu_{33} & 0      & 0      & 0 \\
0      & 0      & 0      & \mu_{44} & 0      & 0\\ 
0      & 0      & 0      & 0      & \mu_{44} & 0\\
0      & 0      & 0      & 0      & 0      & \mu_{66}
\end{pmatrix}
$$
with $\mu_{66}=(\mu_{11}-\mu_{12})/2$. Here
$$ 
\Gc=
\begin{pmatrix}
c_{11}\p_1^2+c_{66}\p_2^2+c_{44}\p_3^2 & \left(c_{11}-c_{66}\right)\p_1\p_2& 0\\
\left(c_{11}-c_{66}\right)\p_1\p_2 & c_{66}\p_1^2+c_{11}\p_2^2+c_{44}\p_3^2 & 0\\
0 & 0 & c_{44}\p_1^2+c_{44}\p_2^2+c_{33}\p_3^2\\
\end{pmatrix}
$$
The eigenvalues $L_i^c(\na)$ of $\Gc(\na)$ in this case are  
\begin{eqnarray*}
L_1^c(\na)=c_{44}\p_1^2+c_{44}\p_2^2+c_{33}\p_3^2\\
L_2^c(\na)=c_{11}\p_1^2+c_{11}\p_2^2+c_{44}\p_3^2\\
L_3^c(\na)=c_{66}\p_1^2+c_{66}\p_2^2+c_{44}\p_3^2
\end{eqnarray*}
and the associated eigenvectors $\D_i^c(\na)$ are:
\begin{eqnarray*}
\D_1^c=\begin{pmatrix}0 \\ 0\\ 1\end{pmatrix}, \quad
\D_2^c=\begin{pmatrix}\p_1\\\p_2\\0\end{pmatrix},\quad
\D_3^c=\begin{pmatrix}\p_2\\-\p_1\\0\end{pmatrix}.
\end{eqnarray*}
Thus $M_1^c=1,\quad\text{ and }\quad M_2^c=M_3^c=\p_1^2+\p_2^2$

\subsection{Medium III}

\par
Finally, we will present the elastodynamic Green function for
another transversely isotropic media with the axis of symmetry 
along ${\bef}_3$ and having $\ul c$ and $\ul v$ of the form:
$$
%{\bf c}=
\begin{pmatrix}
\mu_{11} & \mu_{11}-2\mu_{66} & \mu_{11}-2\mu_{44} & 0      & 0      & 0\\
\mu_{11}-2\mu_{66} & \mu_{11} & \mu_{11}-2\mu_{44} & 0      & 0      & 0\\
\mu_{11}-2\mu_{44} & \mu_{11}-2\mu_{44} & \mu_{11} & 0      & 0      & 0 \\
0      & 0      & 0      & \mu_{44} & 0      & 0\\ 
0      & 0      & 0      & 0      & \mu_{44} & 0\\
0      & 0      & 0      & 0      & 0      & \mu_{66}
\end{pmatrix}
$$
The Christoffel tensor in this case is 
$$
\Gc=
\begin{pmatrix}
c_{11}\p_1^2+c_{66}\p_2^2+c_{44}\p_3^2 & \left(c_{11}-c_{66}\right)\p_1\p_2& \left(c_{11}-c_{44}\right)\p_1\p_3\\
\left(c_{11}-c_{66}\right)\p_1\p_2 & c_{66}\p_1^2+c_{11}\p_2^2+c_{44}\p_3^2 & \left(c_{11}-c_{44}\right)\p_2\p_3\\
\left(c_{11}-c_{44}\right)\p_1\p_3 & \left(c_{11}-c_{44}\right)\p_2\p_3 & c_{44}\p_1^2+c_{44}\p_2^2+c_{11}\p_3^2\\
\end{pmatrix}
$$
Its  eigenvalues $L_i^c(\na)$ are:  
\begin{eqnarray*}
L_1^c(\na)&=&c_{11}\p_1^2+c_{11}\p_2^2+c_{11}\p_3^2=c_{11}\Delta_x\\
L_2^c(\na)&=&c_{66}\p_1^2+c_{66}\p_2^2+c_{44}\p_3^2\\
L_3^c(\na)&=&c_{44}\p_1^2+c_{44}\p_2^2+c_{44}\p_3^2=c_{44}\Delta_x
\end{eqnarray*}
and the eigenvectors $\D_i^c(\na)$ are:
\be 
\D_1^c=\begin{pmatrix}\p_1 \\ \p_2\\\p_3 \end{pmatrix}
,\quad\D_2^c=\begin{pmatrix}\p_2\\-\p_1\\0\end{pmatrix}
,\quad\D_3^c=\begin{pmatrix}-\p_1\p_3\\-\p_2\p_3\\\p_1^2+\p_2^2\end{pmatrix}
\ee
In this case, $M_1^c=\Delta_x\quad M_2^c=\p_1^2+\p_2^2\quad\text{and}\quad M_3^c=(\p_1^2+\p_2^2)\Delta_x$

\subsection{Properties of the Media and Main Assumptions}\label{Assumptions}

\par 
In all anisotropic media discussed above, it holds that 
\begin{itemize}
 \item The Christoffel tensors $\Gc$ and $\Gv$ have the same structure in the sense that 
$$
\D_i^c = \D_i^v,\quad\forall i=1,2,3.
$$
 \item The eigenvalues $L_i^c(\na)$  are homogeneous quadratic forms in the components of the 
argument vector $\na$ \emph{i.e.}
$$
L_i^c[\na]=\dis\sum_j^3a_{i j}^2\frac{\p^2}{\p x_j^2},
$$ 
and therefore  equations  \eqref{eq:PHI} are actually scalar wave equations.
\item In all the concerning cases, the operator  $M_i^c (\na)$ is 
either constant or has a homogeneous quadratic form 
$$
M_i^c=\dis\sum_j^3m_{ij}^2\frac{\p^2}{\p x_j^2}.
$$ 
\end{itemize}

\par 
In addition, we assume that
\begin{itemize}
 \item the eigenvalues of $\Gc$ and $\Gv$ satisfy
$$ L_i^v(\na) = \beta_i  L_i^c(\na).$$
\item and the loss per wave length is small, i.e. 
$$ \beta_i << 1.$$
\end{itemize}

\begin{remark}
 The expression  $ M_3^c=(\p_1^2+\p_2^2)\Delta_x$ will be avoided in the construction of
the Green function by using the expression 
$$ 
\bG=\phi_3\ul I+\E_1(\na)(\phi_1-\phi_3)+\E_2(\na)(\phi_2-\phi_3)
$$ 
for the elastodynamic Green function.
\end{remark}

\begin{remark}\label{RemD}
In general, $\D_i^c$ and $\D_i^v$ are dependant on the parameters  $c_{pq}$ and $v_{pq}$. Consequently,  $\Gc$ and $\Gv$ can not be diagonalized simultaneously. However, in certain restrictive cases where the polarization directions of different wave modes (\emph{i.e.} quasi longitudinal (qP) and quasi shear waves (qSH and qSV)) are independent of the stiffness or viscosity parameters, it is possible to diagonalize both $\Gc$ and $\Gv$ simultaneously.  Moreover, the assumption on the eigenvalues $L_i^v$ and $L_i^c$, implies that for a given  wave mode, the decay rate of its velocity in different directions is uniform, but for different wave modes (qP, qSH and qSV) these decay rates are different. 
\end{remark}

\section{Solution of the Model Wave Problem}\label{WaveProblem}%%%%

Let us now study the scalar wave problems \eqref{eq:PHI}. We consider a
model problem and drop the subscript for brevity in this section
 as well as in the next section. Consider
\be 
\left( L^c[\na]\phi +L^v[\na]\A[\phi] \right)-\dis\rho\frac{\p^2\phi}{\p t^2}
=\delta(t)\delta(\bx).
\label{eq:PHIgeneral}
\ee

\par 
Our assumptions  on the media imply that $L^c$ and $L^v$ have the following form; 
$$ 
L^c[\na]=\dis\sum_{j=1}^3 a_j^2\frac{\p^2}{\p x_j^2}
\quad\text{and}\quad
L^v[\na]=  \beta L^c[\na] = \dis\sum_{j=1}^3 \beta a_j^2\frac{\p^2}{\p x_j^2}
$$
Therefore, the model equation $(\ref{eq:PHIgeneral})$ can be rewritten as:
$$
\dis\sum_{j=1}^3\left(a_j^2\frac{\p^2\phi}{\p x_j^2}+
\beta a_j^2\A\left[\frac{\p^2\phi}{\p x_j^2}\right]\right)
-\dis\rho\frac{\p^2\phi}{\p t^2}=\delta(t)\delta(\x), 
$$

\par 
By a change of variables $x_j=\dis\frac{a_j}{\sqrt\rho}\xi_j$,  we obtain in function  $\tphi(\xi)=\phi(\x)$ the following  transformed equation :
\be 
\Delta_\xi\tphi+\beta\A\left[\Delta_\xi\tphi\right]-\frac{\p^2\tphi}{\p t^2}
=\frac{\sqrt\rho}{a}\delta(t)\delta(\xi).\label{TransfromedEqn}
\ee
where the constant $a=a_1a_2a_3$.

\par 
Now, we apply $\A$ on both sides of the equation \eqnref{TransfromedEqn}, 
and replace the resulting  expression for $\A\left[\Delta_ \xi\tphi\right]$ 
back  into the equation \eqnref{TransfromedEqn}. This yields: 
$$ 
\Delta_\xi\tphi+\beta\A\left[\frac{\p^2\tphi}{\p t^2}\right]
-\beta^2\A^2\left[\Delta_\xi\tphi\right]-\frac{\p^2\tphi}{\p t^2}
=\frac{\sqrt\rho}{a}\delta(\xi)\left\{\delta(t)-\beta\A[\delta(t)]\right\}
$$
Recall that  $\beta << 1$ and  the term in $\beta^2$ is negligible. Therefore, it holds
\be 
\Delta_\xi\tphi+\beta\A\left[\dis\frac{\p^2\tphi}{\p t^2}\right]
-\dis\frac{\p^2\tphi}{\p t^2}
\simeq
\frac{\sqrt\rho}{a}\delta(\xi)\left\{\delta(t)-\beta\A[\delta(t)]\right\}.
\label{PhiApprox}
\ee
Finally, taking temporal Fourier transform 
on both sides of \eqref{PhiApprox}, we obtain the 
corresponding Helmholtz equation: 
\be 
\Delta_\xi\tPhi+\omega^2\left(1-\beta\widehat\A(\omega)\right)\tPhi
=\left(1-\beta\widehat\A(\omega)\right)\frac{\sqrt\rho}{a}\delta(\xi)\label{Helmholtz}
\ee
where $\widetilde\Phi(\xi,\omega)$ 
and $\widehat\A(\omega)$ are the Fourier transforms of
 $\widetilde\phi(\xi,t)$ and the kernel of the convolution operator $\A$ 
respectively.
Let 
$$
\K(\omega)=\dis\sqrt{\omega^2\left(1-\beta\widehat\A(\omega)\right)}.
$$
Then the  solution of the Helmholtz equation \eqnref{Helmholtz} (see for instance \cite{CourantHilbert, Nedelec}) is expressed as 
$$ 
\Phi(\x,\omega)=\sqrt\rho\left(1-\beta\widehat\A(\omega)\right)
\frac{e^{\sqrt{-1}\K(\omega)\tau(\x)}}{4a\pi\tau(\x)}.
$$
where 
$$
\tau(\x)=\sqrt\rho\sqrt{\frac{x_1^2}{a_1^2}+ \frac{x_2^2}{a_2^2}+
 \frac{x_3^2}{a_3^2}}
$$
Using density normalized constants $b_j=\dis\frac{a_j}{\sqrt\rho}$,
 we have  
\be 
\Phi(\x,\omega)
=\left(1-\beta\widehat\A(\omega)\right)
\frac{e^{\sqrt{-1}\K(\omega)\tau(\x)}}{4b\rho\pi\tau(\x)}.
\ee 
where constant $b=b_1b_2b_3$ and
$$
\tau(\x)=\sqrt{\frac{x_1^2}{b_1^2}+ \frac{x_2^2}{b_2^2}+
 \frac{x_3^2}{b_3^2}}
$$

\section{Solution of the Model Potential Problem}\label{section_potential}

\par 
In this section, we find the solution of equation 
\eqref{PotentialProblem}. We once again proceed with a model 
problem. Once the solution is obtained, we will aim to
calculate, its second order derivatives for the evaluation of 
$\D  \otimes \D \psi$.

\subsection{Solution of the Potential Problem}

\par 
Let $\psi(\x,t)$, be the solution of equation \eqref{PotentialProblem}
and  $\Psi(\bx,\omega)$ be its Fourier transform with respect to 
variable $t$. Then $\Psi(\bx, \omega)$ satisfies, 
\be 
M\Psi(\x,\omega)=\Phi(\x,\omega)=\left(1-\beta\widehat\A(\omega)\right)
\frac{e^{\sqrt{-1}\K(\omega)\tau(\x)}}{4b\rho\pi\tau(\x)}.\label{MPsi}
\ee

\par 
When $M$ is constant,  the solution of this equation is directly calculated. 
As $M=(\p_1^2+\p_2^2)\Delta_x$ will not be used in the construction of
 Green function, we are only interested in the case where $M$ is
a homogeneous quadratic form in the component of $\na$ \emph{i.e.}
$$ 
M=\dis\sum_{j=1}^3 m_j^2\frac{\p^2}{\p x_j^2}.
$$
So, the model equation \eqref{MPsi} can be rewritten as:
\be 
\dis\sum_{j=1}^3m_j^2\frac{\p^2\Psi}{\p x_j^2}
=\left(1-\beta\widehat\A(\omega)\right)
\frac{e^{\sqrt{-1}\K(\omega)\tau(\x)}}{4b\rho\pi \tau(\x)}
\quad m_j\neq 0,\quad \forall j\label{ModelPotential}
\ee

\par
By a change of variables $x_j=m_j\eta_j$, equation \eqref{ModelPotential} 
becomes Poisson equation in $\bPsi(\eta,\omega)=\Psi(\x,\omega)$ 
\emph{i.e.}
\be
\Delta_\eta\bPsi=\left(1-\beta\widehat\A(\omega)\right)
\frac{e^{\sqrt{-1}\K(\omega)\overline\tau(\eta)}}{4b\rho\pi\overline\tau(\eta)}
=\bPhi(\eta,\omega)\label{Poisson}
\ee 
where, 
$$
\overline\tau(\eta)=\sqrt{\frac{m^2_1\eta_1^2}{b_1^2}+ 
\frac{m^2_2\eta_2^2}{b_2^2}+ \frac{m^2_3\eta_3^2}{b_3^2}}=\tau(\x)
\quad\text{and}\quad\bPhi(\eta,\omega)=\Phi(\x,\omega)
$$

\par 

Notice that the source 
$\bPhi(\eta,\omega)$ is symmetric with respect to ellipsoid 
$\overline\tau$, \emph{i.e.}
$$
\bPhi(\eta,\omega)=\bPhi(\overline\tau,\omega).
$$
Therefore, the solution $\bPsi$ of the Poisson 
equation \eqref{Poisson}  is the potential field  of a 
uniformly charged ellipsoid due to a charge density $\bPhi(\overline\tau,\omega)$. 
The potential field $\bPsi$ can be calculated with a classical approach
using ellipsoidal coordinates. (See for example 
\cite{Chandrasekhare, kellogg1953foun}
 for the theory of potential problems in ellipsoidal coordinates.)

\par
For the solution of the Poisson equation \eqref{Poisson}
we recall following result from \cite[Ch. 7, Sec.6]{kellogg1953foun}.

\begin{prop}

Let 
$$ 
f(z)=\dis\sum_{j=1}^3\frac{\zeta_j^2}{\left(\alpha_jh\right)^2+z}-1\quad  
\quad\text{and}\quad 
g(z)=\Pi_{j=1}^3\left[\left(\alpha_jh\right)^2+z\right]
$$
and let $Z(h,\zeta)$ be the  largest algebraic root of $f(z)g(z)=0$.
Then the solution of the Poisson equation 
$$\Delta^2 Y(\zeta)=4\pi 
\chi\left(\frac{\zeta_1^2}{\alpha_1^2}+\frac{\zeta_2^2}{\alpha_2^2}
+\frac{\zeta_2^2}{\alpha_1^2}\right)
\quad\zeta\in\RR^3$$
is given by 
$$
Y(\zeta)=2\pi{\alpha_1\alpha_2\alpha_3}
\int_0^\infty\chi(h) I(h,\zeta)dh.
$$
The integrand $I(h,\zeta)$ is defined as 
\begin{equation*} 
I(h,\zeta)=\left|
\begin{array}{lr}
\dis h^2\int_{Z(h,\zeta)}^\infty\frac{1}{\sqrt{g(z)}}dz& Z>0
\\\\
\dis h^2\int_0^\infty\frac{1}{\sqrt{g(z)}}dz & Z<0
\end{array}
\right.
\end{equation*}

\end{prop}

\par 
Hence, the solution of \eqref{Poisson} can be given as  
\begin{equation*} 
\bPsi(\eta,\omega)=\dis \frac{2\pi b}{m}\left(1-\beta\widehat\A(\omega)\right)
\frac{1}{4\pi}
\int_0^\infty\frac{e^{\sqrt{-1}\K(\omega)h}}{4b\rho\pi  h} I(h,\eta)dh
\end{equation*}
or equivalently, 
\begin{equation} 
\Psi(\x,\omega)=\dis\frac{1}{8\rho\pi m}
\left(1-\beta\widehat\A(\omega)\right)\int_0^\infty
\frac{e^{\sqrt{-1}\K(\omega)h}}{h} I(h,\x)dh,\quad m=m_1m_2m_3
\label{Sol_Poisson}
\end{equation}

\par
By a change of variable $s= h^{-2}z$, we can write $I(h,\x)$ as:
\begin{equation} 
I(h,\x)=\left|
\begin{array}{lr}
\dis mh\int_{S(h,\x)}^\infty\frac{1}{\sqrt{G(s)}}ds& h<\tau
\\\\
\dis mh\int_0^\infty\frac{1}{\sqrt{G(s)}}ds & h>\tau
\end{array}
\right.
\end{equation}
with $S(h,\x)=h^{-2}Z(h,\x)$ being the largest algebraic root of the
 equation 
$$F(s)G(s)=0$$ 
where
\be 
\left|
\begin{array}{l}
F(s)=h^2f(h^2s)=\dis\sum_{j=1}^3 \left\{V_j(s)\right\}^{-1}x_j^2-h^2
\\
\\
G(s)=\dis\frac{m^2}{h^6}g(h^2s)=\Pi_{j=1}^3\left\{V_j(s)\right\} 
\\
\\
\text{with}\quad V_j(s)=b_j^2+m_j^2s
\end{array}
\right.
\ee

\begin{remark} 
Note that, $F(s)\equiv 0$ corresponds to a set of confocal ellipsoids 
\be 
s\longmapsto h^2(s)=\sum_{j=1}^3 \left\{V_j(s)\right\}^{-1}x_j^2
\ee 
such that  $\tau(\x)=h(0)$ \emph{i.e.} $S(\tau)=0$. Moreover, $S>0$ if the ellipsoid
$h$ lies inside $\tau$ and $S<0$ if the ellipsoid $h$ lies outside $\tau$.
\end{remark}

\subsection{Derivatives of the Potential field}

\par 
Now we compute the derivatives of the potential $\Psi$.
We note that $I(h,\x)$ is constant with respect to $\x$ 
when $h>\tau.$ So,
\begin{equation*} 
\dis\frac{\p I(h,\x)}{\p x_k} =\left|
\begin{array}{lr}
\dis -mh\frac{\p S(h,\x)}{\p x_k}\frac{1}{\sqrt{G(S(h,\x))}} & h<\tau
\\\\
0 & h>\tau
\end{array}
\right.
\end{equation*}
for $k=1,2,3$ and by consequence,
$$ 
\begin{array}{l}
\dis\frac{\p \Psi}{\p x_k}=-\dis\frac{1}{8\rho\pi m}
\left(1-\beta\widehat\A(\omega)\right)\int_0^\infty\frac{e^{\sqrt{-1}\K(\omega)h}}{h}
\dis\frac{\p I(h,\x)}{\p x_k}dh
\end{array}
$$
or
\be 
\begin{array}{l}
\dis\frac{\p \Psi}{\p x_k}=-\dis\frac{1}{8\rho\pi}
\left(1-\beta\widehat\A(\omega)\right)
\int_0^\tau\left[e^{\sqrt{-1}\K(\omega)h}\right]\frac{\p S(h,\x)}{\p x_k}
\frac{1}{\sqrt{G(S(h,\x))}}dh.
\end{array}\label{FirstDiff}
\ee
Now, we apply $\dis\frac{\p}{\p x_l}$ for $l=1,2,3$ on \eqref{FirstDiff} to obtain the 
second order derivatives of $\Psi$:
\begin{eqnarray*}
-{8\rho\pi}\dis\frac{\p^2 \Psi}{\p x_kx_l}
&=&
\dis\left(1-\beta\widehat\A(\omega)\right)
\dis\frac{\p}{\p x_l}\left[\int_0^\tau\left[e^{\sqrt{-1}\K(\omega)h}\right]
\frac{\p S}{\p x_k}\frac{1}{\sqrt{G(S)}}dh\right]
\\
&=&\dis\left(1-\beta\widehat\A(\omega)\right)\frac{\p\tau}{\p x_l}
\left\{\left[e^{\sqrt{-1}\K(\omega)\tau}\right]
\frac{\p{S(\tau)}}{\p x_k}\frac{1}{\sqrt{G(S(\tau))}}\right\}
\\
&+& \dis\left(1-\beta\widehat\A(\omega)\right)\int_0^\tau
\left[e^{\sqrt{-1}\K(\omega)h}\right]\frac{1}{\sqrt{G(S)}}
\left\{\frac{\p^2S}{\p x_k\p x_l}-\frac{1}{2}\frac{\p S}{\p x_k}
\frac{\p S}{\p x_l}\frac{G'(S)}{G(S)}\right\}dh
\end{eqnarray*}

\par
As $F(S)G(S)=0$ and  $G(s)$ is normally non-zero on $S$,
 therefore by differentiating $F(S)=0$, we obtain 
\cite[eq. (5.21)-(5.23)]{Burridge93}
\be 
\dis\frac{\p S}{\p x_k}= \frac{-2x_k}{V_k(S)F'(S)}\label{Der_S}
\ee 

\be
\frac{\p^2 S}{\p x_kx_l} 
=
\frac{-4x_kx_l}{V_k(S)V_l(S)\left[{F'(S)}\right]^2}
\left\{\frac{F''(S)}{F'(S)}+\frac{m_k^2}{V_k(S)}+
\frac{m_l^2}{V_l(S)}\right\}
- \frac{2\delta_{kl}}{V_k(S)F'(S)}\label{DDer_S}
\ee
where,
\be 
\begin{array}{l}
F'(s)=\dis\sum_{j=1}^3\dis\frac{-m_j^2x_j^2}{V_j^2(s)},
\quad
F''(s)=\sum_{j=1}^3\dis\frac{2m_j^4x_j^2}{V_j^3(s)},
\quad
G'(s)=G(s)\sum_{j=1}^3\dis\frac{m_j^2}{V_j(s)}
\end{array}
\ee
and prime represents a derivative with respect to variable $s$.

\par 
Substituting the values from \eqref{Der_S} and \eqref{DDer_S}, 
the second order derivative of $\Psi$ becomes

\be 
{4\rho\pi}\dis\frac{\p^2 \Psi}{\p x_kx_l}
 =\left|
\begin{array}{l}
\dis\frac{-x_kx_l\dis\left(1-\beta\widehat\A(\omega)\right)}{aa_k^2a_l^2F'(0)}
\left\{\frac{e^{\sqrt{-1}\K(\omega)\tau}}{\tau}\right\}
\\
+\dis\left(1-\beta\widehat\A(\omega)\right)\int_0^\tau
\left[e^{\sqrt{-1}\K(\omega)h}\right]\frac{1}{F'(S)\sqrt{G(S)}}\times
\\
\dis\left[\frac{2x_kx_l}{V_k(S)V_l(S){F'(S)}}
\left\{\frac{F''(S)}{F'(S)}+\frac{m_k^2}{V_k(S)}+\frac{m_l^2}{V_l(S)}
+\frac{1}{2}\frac{G'(S)}{G(S)}\right\}
+ \dis\frac{\delta_{kl}}{V_k(S)}\right]dh
\end{array}
\right.
\label{SecDerPotential}
\ee 
\begin{remark}
If for some $i\in\{1,2,3\}$, $m_i\to 0$ one semi axis of the ellipsoid $\tau$
 tends to infinity but no singularity occurs. Therefore the results of this 
section are still valid in this case.
\end{remark}

\section{Elastodynamic Green Function}

\begin{table}[b]
\centering
\begin{tabular}{|c|ccc|ccc|c|} 
\hline
Medium &$b_1$ &$b_2$&$b_3$ &$m_1$ &$m_2$& $m_3$& $M_i$ \\
\hline
\multirow{3}{*}{I} & $c_1$ & $c_6$ & $c_5$ & 1 & 0 & 0 & $M_1$ \\
                   & $c_6$ & $c_2$ & $c_4$ & 0 & 1 & 0 & $M_2$ \\
                   & $c_5$ & $c_4$ & $c_3$ & 0 & 0 & 1 & $M_3$ \\
\hline
\multirow{3}{*}{II}& $c_4$ & $c_4$ & $c_3$ & 0 & 0 & 1 & $M_1$ \\
                   & $c_1$ & $c_1$ & $c_4$ & 1 & 1 & 0 & $M_2$ \\
                   & $c_6$ & $c_6$ & $c_4$ & * & * & * & $M_3$ \\
\hline
\multirow{3}{*}{III} & $c_1$ & $c_1$ & $c_1$ & 1 & 1 & 1 & $M_1$ \\
                   & $c_6$ & $c_6$ & $c_4$ & 1 & 1 & 0 & $M_2$ \\
                   & $c_4$ & $c_4$ & $c_4$ & * & * & * & $M_3$ \\
\hline
\end{tabular}
\vspace{1em}
\caption{Values of $b_i$ and $m_i$ for different media. Here $*$ represents 
a value which is not used for reconstructing Green function.}
\label{Table}
\end{table}

\par 
In this section we present the expressions for the elastodynamic Green
functions for the media presented in section \ref{sec:Media}. 
Throughout this section  $c_{p}=\dis\sqrt{\frac{c_{pp}}{\rho}}$ 
with $p\in\{1,2,\cdots,6\}$. We recall that
$
\K_i(\omega)=\dis\sqrt{\omega^2\left(1-\beta_i\widehat\A(\omega)\right)}.
$
\subsection{Medium I}
All the eigenvectors of $\G$ are constants in this case \emph{i.e.}
 $\D_i=\bef_i$, therefore $M_i=1$ and $\E_i=\bef_i\otimes\bef_i$. 
If $\widehat\bG$ is the Fourier transform  of the 
viscoelastic Green function $\bG$ with respect to variable $t$, then: 
\be 
\widehat\bG=\dis\sum_{i=1}^3 \Phi_i(x,\omega)\bef_i\otimes\bef_i=
\frac{1}{4\pi \rho}\sum_{i=1}^3 
\left[
\frac{c_{i+3}\left(1-\beta_i\widehat\A(\omega)\right)}{c_ic_4c_5c_6\tau_i}
\exp(\sqrt{-1}\K_i(\omega)\tau_i)\right]\bef_i\varotimes\bef_i
\ee 
where 
$$
\tau_1=\sqrt{\frac{x_1^2}{c_1^2}+\frac{x_2^2}{c_6^2}+\frac{x_3^2}{c_5^2}}
,\quad
\tau_2=\sqrt{\frac{x_1^2}{c_6^2}+\frac{x_2^2}{c_2^2}+\frac{x_3^2}{c_4^2}}
,\quad
\tau_3=\sqrt{\frac{x_1^2}{c_5^2}+\frac{x_2^2}{c_4^2}+\frac{x_3^2}{c_3^2}}
,\quad
$$

\subsection{Medium II}

According to section \ref{WaveProblem}, the functions 
$\Phi_i$ have following expressions: 
\begin{eqnarray*} 
\Phi_1(\x,\omega)
=\dis\left(1-\beta_1\widehat\A(\omega)\right)
\frac{e^{\sqrt{-1}\K_1(\omega)\tau_1(\x)}}{4c_4^2c_3\rho\pi\tau_1(\x)}
\\
\Phi_2(\x,\omega)
=\dis\left(1-\beta_2\widehat\A(\omega)\right)
\frac{e^{\sqrt{-1}\K_2(\omega)\tau_2(\x)}}{4c_1^2c_4\rho\pi\tau_2(\x)}
\\
\Phi_3(\x,\omega)
=\dis\left(1-\beta_3\widehat\A(\omega)\right)
\frac{e^{\sqrt{-1}\K_3(\omega)\tau_3(\x)}}{4c_6^2c_4\rho\pi\tau_3(\x)}
\end{eqnarray*}
where 
$$
\tau_1(\x)=\sqrt{\frac{x_1^2}{c_4^2}+ \frac{x_2^2}{c_4^2}+
 \frac{x_3^2}{c_3^2}},
\quad
\tau_2(\x)=\sqrt{\frac{x_1^2}{c_1^2}+ \frac{x_2^2}{c_1^2}+
 \frac{x_3^2}{c_4^2}},
\quad
\tau_3(\x)=\sqrt{\frac{x_1^2}{c_6^2}+ \frac{x_2^2}{c_6^2}+
 \frac{x_3^2}{c_4^2}}
$$

To calculate Green function, we use the expression
$$
\widehat\bG=\Phi_3\ul I+\dis\D_1\otimes\D_1M_1^{-1}\left(\Phi_1-\Phi_3\right)
+\dis\D_2\otimes\D_2M_2^{-1}\left(\Phi_2-\Phi_3\right).
$$
$\D_1=\bef_3$ and $M_1=1$, yield  
$$
\D_1\otimes\D_1M_1^{-1}\left(\Phi_1-\Phi_3\right)
=\left(\Phi_1-\Phi_3\right)\bef_3\otimes\bef_3
$$
To compute 
$\D_2\otimes\D_2M_2^{-1}\left(\Phi_2-\Phi_3\right)$, suppose 
$$
\Psi_2=M_2^{-1}\Phi_2\quad\text{and}\quad\Psi_3=M_2^{-1}\Phi_3
$$
and notice that $m_1=m_2=1$ and $m_3=0$. Moreover for $\Phi_2$ and $\Phi_3$, $b_1=b_2$. (See Table \ref{Table}). Thus, we have
\begin{equation*} 
\begin{array}{l}
\dis\frac{4\rho\pi}{\left(1-\beta_2\widehat\A(\omega)\right)}
\dis\frac{\p^2 \Psi_2}{\p x_kx_l}
 =
{\hR_k\hR_l}\left\{\frac{e^{\sqrt{-1}\K_2(\omega)\tau_2}}{c_1^2c_4\tau_2}\right\}
-\dis\frac{1}{c_4R^2}(\delta_{kl}- 2\hR_k\hR_l) 
\int_0^{\tau_2}\left[e^{\sqrt{-1}\K_2(\omega)h}\right]dh
\\
\\
\dis\frac{4\rho\pi}{\left(1-\beta_3\widehat\A(\omega)\right)}
\dis\frac{\p^2 \Psi_3}{\p x_kx_l}
 =
{\hR_k\hR_l}\left\{\frac{e^{\sqrt{-1}\K_3(\omega)\tau_3}}{c_6^2c_4\tau_3}\right\}
-\dis\frac{1}{c_4R^2}(\delta_{kl}- 2\hR_k\hR_l) 
\int_0^{\tau_3}\left[e^{\sqrt{-1}\K_3(\omega)h}\right]dh
\end{array}
\end{equation*}
where $\hR_k=\dis\frac{x_k}{R}$ for $k=1,2$. 
See Appendix \ref{appendix_II} for the derivation of this result.

\par 
By using the second derivatives of $\Psi_2$ and $\Psi_3$
and the  expression
$$
%\begin{array}{l}
\D_2\otimes\D_2M_2^{-1}\left(\Phi_2-\Phi_3\right)
%\\\\=
%\dis\left(\p_1\bef_1+\p_2\bef_2\right)
%\otimes\left(\p_1\bef_1+\p_2\bef_2\right)\left(\Psi_2-\Psi_3\right)
=\dis\sum_{k,l=1}^2\p_k\p_l\left(\Psi_2-\Psi_3\right)\bef_k\otimes\bef_l
%\end{array}
$$
we finally arrive at
$$
\begin{array}{l}
\widehat\bG
=\dis\left(1-\beta_3\widehat\A(\omega)\right)
\frac{e^{\sqrt{-1}\K_3(\omega)\tau_3(\x)}}{4c_6^2c_4\rho\pi\tau_3(\x)}
\ul J
+\dis\left(1-\beta_1\widehat\A(\omega)\right)
\frac{e^{\sqrt{-1}\K_1(\omega)\tau_1(\x)}}{4c_4^2c_3\rho\pi\tau_1(\x)}
\bef_3\otimes\bef_3
\\
\\
+\left[\dis\left(1-\beta_2\widehat\A(\omega)\right)
\frac{e^{\sqrt{-1}\K_2(\omega)\tau_2(\x)}}{4c_1^2c_4\rho\pi\tau_2(\x)}
-\dis\left(1-\beta_3\widehat\A(\omega)\right)
\frac{e^{\sqrt{-1}\K_3(\omega)\tau_3(\x)}}{4c_6^2c_4\rho\pi\tau_3(\x)}
\right]
\hR\otimes\hR
\\
\\
-\dis\frac{1}{4\rho\pi c_4R^2}(\ul J- 2\hR\otimes\hR)\times
\\
\\
\left[ 
\left(1-\beta_2\widehat\A(\omega)\right)\dis\int_0^{\tau_2}
\left[e^{\sqrt{-1}\K_2(\omega)h}\right]dh
-\left(1-\beta_3\widehat\A(\omega)\right)\dis\int_0^{\tau_3}
\left[e^{\sqrt{-1}\K_3(\omega)h}\right]dh
\right]
\end{array}
$$
Or equivalently,  
$$
\begin{array}{l}
\widehat\bG
=\Phi_1\bef_3\otimes\bef_3
+\Phi_2\hR\otimes\hR
+\Phi_3(\ul J-\hR\otimes\hR)
\\
\\
-\dis\frac{1}{R^2}
\left[c_1^2\int_0^{\tau_2}h\Phi_2(h,\omega)dh
-
c_6^2\int_0^{\tau_3}h\Phi_3(h,\omega)dh
\right]
(\ul J- 2\hR\otimes\hR)
\end{array}
$$
Here $\ul J=\ul I-\bef_3\otimes\bef_3$ and 
$\hR=\hR_1\bef_1+\hR_2\bef_2$

\subsection{Medium III}

\par 
The solutions of the wave equation $\Phi_i$ in this case are
\begin{eqnarray*} 
\Phi_1(\x,\omega)
=\dis\left(1-\beta_1\widehat\A(\omega)\right)
\frac{e^{\sqrt{-1}\K_1(\omega)\tau_1(\x)}}{4c_1^3\rho\pi\tau_1(\x)}
\\
\Phi_2(\x,\omega)
=\dis\left(1-\beta_2\widehat\A(\omega)\right)
\frac{e^{\sqrt{-1}\K_2(\omega)\tau_2(\x)}}{4c_6^2c_4\rho\pi\tau_2(\x)}
\\
\Phi_3(\x,\omega)
=\dis\left(1-\beta_3\widehat\A(\omega)\right)
\frac{e^{\sqrt{-1}\K_3(\omega)\tau_3(\x)}}{4c_4^3\rho\pi\tau_3(\x)}
\end{eqnarray*}
where 
$$
\tau_1(\x)=\frac{1}{c_1}\sqrt{x_1^2+x_2^2+x_3^2}=\frac{r}{c_1},
\quad
\tau_2(\x)=\sqrt{\frac{x_1^2}{c_6^2}+ \frac{x_2^2}{c_6^2}+
 \frac{x_3^2}{c_4^2}},
\quad
\tau_3(\x)=\frac{r}{c_4}
$$

\par 
To calculate Green function, we once again use the expression
$$
\widehat\bG=\Phi_3\ul I+\dis\D_1\otimes\D_1M_1^{-1}\left(\Phi_1-\Phi_3\right)
+\dis\D_2\otimes\D_2M_2^{-1}\left(\Phi_2-\Phi_3\right).
$$
Suppose $\Psi_1=M_1^{-1}\Phi_1\quad\text{and}\quad\Psi_3=M_1^{-1}\Phi_3$.
Notice that $m_1=m_2=m_3=1$ for $M_1$ and  $b_1=b_2=b_3$ for $\Phi_1$  as well as $\Phi_3$ (see Table \ref{Table}). Thus,
\begin{equation*} 
\begin{array}{l}
\dis\frac{4\rho\pi}{\left(1-\beta_1\widehat\A(\omega)\right)}
\dis\frac{\p^2 \Psi_1}{\p x_kx_l}
 ={\hr_k\hr_l}
\left\{\frac{e^{\sqrt{-1}\K_1(\omega)\tau_1}}{c_1^3\tau_1}\right\}
-\dis\frac{1}{r^3}(\delta_{kl}- 3\hr_i\hr_j)
 \int_0^{\tau_1}\left[he^{\sqrt{-1}\K_1(\omega)h}\right]dh
\\
\\
\dis\frac{4\rho\pi}{\left(1-\beta_3\widehat\A(\omega)\right)}
\dis\frac{\p^2 \Psi_3}{\p x_kx_l}
 ={\hr_k\hr_l}
\left\{\frac{e^{\sqrt{-1}\K_3(\omega)\tau_1}}{c_1^3\tau_3}\right\}
-\dis\frac{1}{r^3}(\delta_{kl}- 3\hr_i\hr_j)
 \int_0^{\tau_3}\left[he^{\sqrt{-1}\K_3(\omega)h}\right]dh
\end{array}
\end{equation*}
See Appendix \ref{appendix_I} for the derivation of this result. 
It yields
$$
\begin{array}{l}
\D_1\otimes\D_1M_1^{-1}\left(\Phi_1-\Phi_3\right)
\\
\\
=\dis\frac{1}{4\rho\pi}\left[\dis\left(1-\beta_1\widehat\A(\omega)\right)
\frac{e^{\sqrt{-1}\K_1(\omega)\tau_1(\x)}}{c_1^3\tau_1(\x)}
+
\dis\left(1-\beta_3\widehat\A(\omega)\right)
\frac{e^{\sqrt{-1}\K_3(\omega)\tau_3(\x)}}{c_4^3\tau_3(\x)}
\right]\hr\otimes\hr
\\
\\
-
\left[
\dis\left(1-\beta_1\widehat\A(\omega)\right)
\int_0^{\tau_1}\left[he^{\sqrt{-1}\K_1(\omega)h}\right]dh
-\dis\left(1-\beta_3\widehat\A(\omega)\right)
\int_0^{\tau_3}\left[he^{\sqrt{-1}\K_3(\omega)h}\right]dh
\right]\times
\\
\\
\dis\frac{1}{4\rho\pi r^3}\left(\ul I-3\hr\otimes\hr\right)
\\
\\
=\dis\left[\Phi_1(x,\omega)-\Phi_3(x,\omega)\right]\hr\otimes\hr
-
\frac{1}{r^3}\left[\dis\int_0^{\tau_1}h^2\Phi_1(h,\omega)dh
-\dis\int_0^{\tau_3}h^2\Phi_3(h,\omega)dh
\right]\dis\left(\ul I-3\hr\otimes\hr\right)
\end{array}.
$$
where $\hr=\hr_1\bef_1+\hr_2\bef_2+\hr_3\bef_3$ with 
$\hr_i=\dis\frac{x_i}{r}$ for all $i=1,2,3$.

\par 
To compute, $\D_2\otimes\D_2M_2^{-1}\left(\Phi_2-\Phi_3\right)$, suppose 
$\Psi_2=M_2^{-1}\Phi_2$ and $\Psi_4=M_2^{-1}\Phi_3$.
By using formula \eqref{CaseII} with $m_1=m_2=1$ and $m_3=0$, we obtain:
\begin{equation*} 
\begin{array}{l}
\dis\frac{4\rho\pi}{\left(1-\beta_2\widehat\A(\omega)\right)}
\dis\frac{\p^2 \Psi_2}{\p x_kx_l}
 =
{\hR_k\hR_l}\left\{\frac{e^{\sqrt{-1}\K_2(\omega)\tau_2}}{c_6^2c_4\tau_2}\right\}
-\dis\frac{1}{c_4R^2}(\delta_{kl}- 2\hR_k\hR_l) 
\int_0^{\tau_2}\left[e^{\sqrt{-1}\K_2(\omega)h}\right]dh
\\
\\
\dis\frac{4\rho\pi}{\left(1-\beta_3\widehat\A(\omega)\right)}
\dis\frac{\p^2 \Psi_4}{\p x_kx_l}
 =
{\hR_k\hR_l}\left\{\frac{e^{\sqrt{-1}\K_3(\omega)\tau_3}}{c_4^3\tau_3}\right\}
-\dis\frac{1}{c_4R^2}(\delta_{kl}- 2\hR_k\hR_l) 
\int_0^{\tau_3}\left[e^{\sqrt{-1}\K_3(\omega)h}\right]dh
\end{array}
\end{equation*}
with $\hR_k=\dis\frac{x_k}{R}$ and $k, l\in\{1,2\}$. This allows us to write
$$
\begin{array}{l}
\D_2\otimes\D_2M_2^{-1}\left(\Phi_2-\Phi_3\right)
\\
\\
=\dis\frac{1}{4\rho\pi}\left[\dis\left(1-\beta_2\widehat\A(\omega)\right)
\frac{e^{\sqrt{-1}\K_2(\omega)\tau_2(\x)}}{c_1^3\tau_2(\x)}
+
\dis\left(1-\beta_3\widehat\A(\omega)\right)
\frac{e^{\sqrt{-1}\K_3(\omega)\tau_3(\x)}}{c_4^3\tau_3(\x)}
\right]\times
\\
\\ 
\left(\hR_2^2\bef_1\otimes\bef_1-
\hR_1\hR_2[\bef_1\otimes\bef_2+\bef_2\otimes\bef_1]+\hR_1^2
\bef_2\otimes\bef_2
\right)
\\
\\
-
\dis\frac{1}{4c_4\rho \pi R^2}\left[
\dis\left(1-\beta_2\widehat\A(\omega)\right)
\int_0^{\tau_2}\left[e^{\sqrt{-1}\K_2(\omega)h}\right]dh
-\dis\left(1-\beta_3\widehat\A(\omega)\right)
\int_0^{\tau_3}\left[e^{\sqrt{-1}\K_3(\omega)h}\right]dh
\right]\times
\\
\\
\left((1-2\hR_2^2)\bef_1\otimes\bef_1-
2\hR_1\hR_2[\bef_1\otimes\bef_2+\bef_2\otimes\bef_1]+(1-2\hR_1^2)
\bef_2\otimes\bef_2\right)
\\
\\
=\dis\left[\Phi_2(x,\omega)-\Phi_3(x,\omega)\right]\hR^\perp\otimes\hR^\perp
\\
\\
-
\dis\frac{1}{R^2}\left[\dis c_6^2\int_0^{\tau_2}h\Phi_2(h,\omega)dh
-\dis c_4^2\int_0^{\tau_3}h\Phi_3(h,\omega)dh
\right]\dis\left(\ul J-2\hR^\perp\otimes\hR^\perp\right)
\end{array}.
$$
where $\hR^\perp=\hR_2\bef_1-\hR_1\bef_2$ and $\ul J=\ul I-\bef_3\otimes\bef_3$.

\par
Finally, we arrive at
$$
\begin{array}{l}
\widehat\bG
=\Phi_1\hr\otimes\hr
+\Phi_2\hR^\perp\otimes\hR^\perp
+\Phi_3(\ul I-\hr\otimes\hr-\hR^\perp\otimes\hR^\perp)
\\
\\
-
\dis\frac{1}{r^3}\left[\dis\int_0^{\tau_1}h^2\Phi_1(h,\omega)dh
-\dis\int_0^{\tau_3}h^2\Phi_3(h,\omega)dh
\right]\dis\left(\ul I-3\hr\otimes\hr\right)
\\
\\
-\dis\frac{1}{R^2}
\left[c_1^2\int_0^{\tau_2}h\Phi_2(h,\omega)dh
-
c_6^2\int_0^{\tau_3}h\Phi_3(h,\omega)dh
\right]
(\ul J-2\hR^\perp\otimes\hR^\perp)
\end{array}
$$

\subsection{Isotropic Medium}

\par 
When $c_{66}=c_{44}$, medium III becomes isotropic. In this case 
$$
\Phi_2(\x,\omega)=\Phi_3(\x,\omega),
\quad
\beta_1=\beta_2,
\quad
\tau_1(\bx)=\frac{r}{c_1},\quad\text{ and }\quad
\tau_2(\bx)=\frac{r}{c_4}=\tau_3(\bx)
$$
Thus, the Green function in an isotropic medium with 
independent elastic parameters $c_{11}$ and $c_{44}$ can be given in frequency 
domain as:
$$
\begin{array}{l}
\widehat\bG
=\Phi_2\ul I + \D_1\otimes\D_1 M_1^{-1}(\Phi_1-\Phi_2)
\\
\\
=
\Phi_1\hr\otimes\hr
+\Phi_2(\ul I-\hr\otimes\hr)
-
\dis\frac{1}{r^3}\left[\dis\int_0^{\frac{r}{c_1}}h^2\Phi_1(h,\omega)dh
-\dis\int_0^{\frac{r}{c_4}}h^2\Phi_2(h,\omega)dh
\right]\dis\left(\ul I-3\hr\otimes\hr\right)
\end{array}
$$
where $\Phi_1$ and $\Phi_2$ are the same as in medium III. This expression of the Green function has already been reported in a previous work \cite{ELA11}.

\section*{Acknowledgement}

We would like to thank Prof. Habib Ammari (\'Ecole Normale Sup\'erieure-Paris) for his continuous support and encouragement.  We would also like to thank Meisam Sharify (\'Ecole Polytechnique-Paris) for his help and remarks. E. Bretin was supported by the foundation \emph{Digitio}  of France, in terms of 
a Post-doctoral Fellowship. A. Wahab was supported by \emph{Higher Education Commission} of Pakistan in terms of a doctoral fellowship.

\appendix

\section{Decomposition of the Green function} \label{appendix_phi}

\par 
Consider the elastic equation satisfied by $\bG$:
\be 
\left( \Gc(\na) \bG(\bx,t) +  \Gv(\na) \A[\bG](\bx,t)  \right)-\rho\dis\frac{\p^2\bG(\bx,t)}{\p t^2}
=\delta(t)\delta(\bx){\ul I}.\label{GreenEq}  
\ee 
If   $\bG$ is given in the form
\be 
\bG=\dis\sum_{i=1}^3\E_i(\na)\phi_i \label{GreenGamma}
\ee
Then substituting \eqref{GreenGamma} in \eqref{GreenEq} yields:
\begin{eqnarray*}
\delta(t)\delta(\bx){\ul I}
&=&
\left( \Gc(\na)\bG(\bx,t) +  \Gv(\na) A[\bG](\bx,t) \right)
-\rho\dis\frac{\p^2\bG(\bx,t)}{\p t^2}
\\   
&=&  
\dis\sum_{i,j=1}^3 \left( L^c_j(\na)\phi_i +  L^v_j(\na) \A[\phi_i]\right) \E_j(\na)\E_i(\na) 
-\rho\sum_{i=1}^3 \E_i(\na)\frac{\p^2\phi_i(\bx,t)}{\p t^2}
\end{eqnarray*}
By definition $\E_i(\na)$ is a projection operator which satisfies
$$
\E_i(\na)\E_j(\na)=\delta_{ij}E_j(\na)
$$
Consequently, we can have
\begin{eqnarray*}
\delta(t)\delta(\bx){\ul I}   
&=&\dis\sum_{i,j=1}^3 \E_j(\na)\delta_{ij} \rho^{-1}  \left(L^c_j(\na)\phi_i + L^v_j(\na)\A[\phi_i] \right)  
-\rho\sum_{i=1}^3 \E_i(\na)\frac{\p^2\phi_i(\bx,t)}{\p t^2}
\\
&=&  \dis\sum_{i=1}^3\E_i(\na)\left(\left(L^c_i(\na)\phi_i + L^v_i(\na)\A[\phi_i] \right) 
-\rho\frac{\p^2\phi_i(\bx,t)}{\p t^2}\right).
\end{eqnarray*}
Moreover $\ul I=\dis\sum_{i=1}^3\E_i(\na)$, therefore
\begin{eqnarray*}
\dis\sum_{i=1}^3\E_i(\na)\left(\left(L^c_i(\na)\phi_i + L^v_i(\na)\A[\phi_i] \right) 
-\rho\frac{\p^2\phi_i(\bx,t)}{\p t^2}-\delta(t)\delta(\bx)
\right)=0
\end{eqnarray*}
Finally, remark that $\bG$ we can express  in the form \eqref{eq:Green} if the functions $\phi_i$ satisfy equation \eqref{eq:PHI}.

\section{Derivative of Potential: Case I}\label{appendix_I}

\par 

If $b_1=b_2=b_3$ and $m_1=m_2=m_3$, we have
\be
\left|\dis
\begin{array}{l}
V_1(s)=V_2(s)=V_3(s)=b_1^2+m_1^2s
\\
\\
F(s)=\dis\sum_{j=1}^3 \frac{x_j^2}{V_1(s)}-h^2
=\frac{r^2}{V_1(s)}-h^2
\\
\\
F'(s)=\dis\sum_{j=1}^3 \frac{-m_1^2x_j^2}{V_1^2(s)}
=\frac{-m_1^2r^2}{V_1^2(s)}\quad\text{and}\quad
F'(0)=\frac{-m_1^2r^2}{b_1^4}
\\
\\
F''(s)=\dis\sum_{j=1}^3 \frac{2m_1^4x_j^2}{V_1^3(s)}
=\frac{2m_1^4r^2}{V_1^3(s)}
\\
\\
G(s)=\dis\left(V_1(s)\right)^3\quad\text{and}\quad
G'(s)=G(s)\dis\frac{3m_1^2}{V_1(s)}
\end{array}
\right.
\label{Case1Der}
\ee
with $r=\sqrt{x_1^2+x_2^2+x_3^2}$. When $F(S)=0$, we have
\be 
\left|
\begin{array}{l}
V_1(S)=\dis\frac{r^2}{h^2},
\\
\\
\dis\left[\frac{1}{V_k(S)V_l(S){F'(S)}}\right]
=\frac{-1}{m_1^2r^2}
\quad\text{and}\quad
\dis\frac{1}{F'(S)\sqrt{G(S)}}=\dis\frac{-1}{m_1^2rh},
\\
\\
\left\{\dis\frac{F''(S)}{F'(S)}+\frac{m_k^2}{V_k(S)}+\frac{m_l^2}{V_l(S)}
+\frac{1}{2}\frac{G'(S)}{G(S)}\right\}
=\dis\frac{3}{2}\frac{m_1^2}{V_1(S)}=\dis\frac{3}{2}\frac{m_1^2h^2}{r^2}
\end{array}
\right.
\label{Case1DerAt-S}
\ee
Substituting \eqref{Case1Der} and \eqref{Case1DerAt-S} in 
\eqref{SecDerPotential} we finally arrive at:
\begin{equation} 
\begin{array}{l}
\dis\frac{4\rho m_1^2\pi}{\left(1-\beta\widehat\A(\omega)\right)}
\dis\frac{\p^2 \Psi}{\p x_kx_l}
 ={\hr_k\hr_l}
\left\{\frac{e^{\sqrt{-1}\K(\omega)\tau}}{b\tau}\right\}
-\dis\frac{1}{r^3}(\delta_{kl}- 3\hr_i\hr_j)
 \int_0^\tau\left[he^{\sqrt{-1}\K(\omega)h}\right]dh
\end{array}
\end{equation}
where $\hr_j=\dis\frac{x_j}{r}$ for all $j=1,2,3.$

\section{Derivative of Potential: Case II}\label{appendix_II}

If $b_1=b_2$,  $m_1=m_2$ and $m_3=0$, we have
\be
\left|\dis
\begin{array}{l}
V_1(s)=V_2(s)=b_1^2+m_1^2s\quad\text{and} V_3(s)= b_3^2
\\
\\
F'(s)=\dis\sum_{j=1}^2 \frac{-m_1^2x_j^2}{V_1^2(s)}
=\frac{-m_1^2R^2}{V_1^2(s)}\quad\text{and}\quad
F'(0)=\frac{-m_1^2R^2}{b_1^4}
\\
\\
F''(s)=\dis\sum_{j=1}^2 \frac{2m_1^4x_j^2}{V_1^3(s)}
=\frac{2m_1^4R^2}{V_1^3(s)}
\\
\\
G(s)=\dis b_3^2\left(V_1(s)\right)^2\quad\text{and}\quad
G'(s)=G(s)\dis\frac{2m_1^2}{V_1(s)}
\end{array}
\right.
\label{Case2Der}
\ee
with $R=\sqrt{x_1^2+x_2^2}$. For all $l,k\in\{1,2\}$, we have
\be 
\left|
\begin{array}{l}
\dis\left[\frac{1}{V_k(S)V_l(S){F'(S)}}\right]
=\frac{-1}{m_1^2R^2}
\quad\text{and}\quad
\\
\\
\dis\frac{1}{F'(S)\sqrt{G(S)}}=\dis\frac{-V(S)}{m_1^2b_3R^2},
\\
\\
\left\{\dis\frac{F''(S)}{F'(S)}+\frac{m_k^2}{V_k(S)}+\frac{m_l^2}{V_l(S)}
+\frac{1}{2}\frac{G'(S)}{G(S)}\right\}
=\dis\frac{m_1^2}{V_1(S)}
\end{array}
\right.
\label{Case2DerAt-S}
\ee
Substituting \eqref{Case2Der} and \eqref{Case2DerAt-S} in 
\eqref{SecDerPotential} and simple calculations, we finally arrive at:
\begin{equation} 
\begin{array}{l}
\dis\frac{4\rho m_1^2\pi}{\left(1-\beta\widehat\A(\omega)\right)}
\dis\frac{\p^2 \Psi}{\p x_kx_l}
 =
{\hR_k\hR_l}\left\{\frac{e^{\sqrt{-1}\K(\omega)\tau}}{b\tau}\right\}
-\dis\frac{1}{b_3R^2}(\delta_{kl}- 2\hR_k\hR_l) 
\int_0^\tau\left[e^{\sqrt{-1}\K(\omega)h}\right]dh
\end{array}\label{CaseII}
\end{equation}
where $\hR_k=\dis\frac{x_k}{R}$ for $k=1,2$.

\bibliographystyle{alpha}

\begin{thebibliography}{99}

\bibitem{Achenbach} J. D. Achenbach, 
\textit{Wave Propagation in Elastic Solids},
North-Holland Publishing Company, Amsterdam, 1973.

\bibitem{AkiRich} K. Aki, P. G. Richards, 
\textit{Quantitative Seismology}, 
\textbf{Vol 1}, 
W.H. Freeman and Co., San Francisco ,  1980.

\bibitem{AlekseevRybak} V. N. Alekseev, S. A. Rybak,
\textit{Equations of state for viscoelastic biological media},
Acoustical Physics, \textbf{48(5)}: (2002), 511--517.


\bibitem{ShapeOpt} G. Allaire, 
\textit{Shape Optimization by the Homogenization Method},
Applied Mathematical Sciences, 
\textbf{146}, Springer-Verlag, New York, 2002. 

\bibitem{AmmariBound} C. Alves,  H. Ammari,
\textit{Boundary integral formulae for the reconstruction of 
imperfections of small diameter in an elastic medium},
SIAM J. on Applied Mathematics,
\textbf{62(1)}: (2001), 94--106.

\bibitem{AmmariBioMed} H. Ammari,
\textit{An introduction to Mathematics of Emerging Biomedical Imaging},
Mathematics \& Applications,
\textbf{ Vol. (62)}, Springer-Verlag, Berlin, 2008.

\bibitem{AmmariDirEast} H. Ammari, P. Calmon, E. Iakovleva,
\textit{Direct elastic imaging of a small inclusion},
SIAM J. Imaging Sci.,
\textbf{1}: (2008), 169--187.

\bibitem{AmmariSepration} H. Ammari, P. Garapon, F. Jouve,
\textit{Separation of scales in elasticity imaging: A numerical study},
Journal of Computational Mathematics,
\textbf{28(3)}: (2010), 354--370.

\bibitem{AmmariOpti} H. Ammari, P. Garapon, F. Jouve, H. Kang, M. Lim,
\textit{A new optimal control approach for the reconstruction of extended inclusions}, Preprint.

\bibitem{AmmariMRElasto} H. Ammari, P. Garapon, H. Kang, H. Lee,
\textit{A method of biological tissues elasticity reconstruction
using magnetic resonance elastography measurements},
Quarterly of Applied Mathematics, 
\textbf{66(1)}:(2008), 139--176.

\bibitem{AmmariEffVis} H. Ammari, P. Garapon, H. Kang, H. Lee,
\textit{ Effective viscosity properties of dilute suspensions 
of arbitrarily shaped particles},
Preprint.

\bibitem{AmmariTransElas} H. Ammari, L. Guadarrama-Bustos,  H. Kang, H. Lee,
\textit{Transient elasticity imaging and time reversal},
Preprint.

\bibitem{AmmariExp}  H. Ammari, H. Kang,
\textit{Expansion methods}, Handbook of Mathematical Methods in Imaging,
Springer-Verlag, New York, 2011.

\bibitem{AmmariRSI} H. Ammari, H. Kang,
\textit{Reconstruction of Small Inhomogeneities from Boundary
 Measurements}, 
Lecture Notes in Mathematics, \textbf{Vol. 1846},
 Springer-Verlag, Berlin, 2004.

\bibitem{AmmariPol} H. Ammari, H. Kang,
\textit{Polarization and Moment Tensors}, Applied Mathematical Sciences,
\textbf{162}, Springer-Verlag, New York, 2008.

\bibitem{Bercoff} J. Bercoff, M. Tanter, M. Muller, M. Fink,
\textit{The role of viscosity in the impulse diffraction
 field of elastic waves induced by the acoustic radiation force},
IEEE Transactions on  Ultrasonics, Ferroelectrics and Frequncy Control,
\textbf{51(11)}: (2004), 1523--1536.

\bibitem{BinSingh} A. Ben-Menahem, S. J. Singh,
\textit{Seismic waves and sources},
Springer-Verlag, 1981.

\bibitem{ELA11} E. Bretin, L. Guadarrama Bustos, A. Wahab,
\textit{On the Green function in visco-elastic media obeying 
a frequency power-law},
Mathematical Methods in the Applied Sciences,\textbf{}(2011),
 DOI: 10.1002/mma.1404.


\bibitem{Burridge93} R. Burridge, P. Chadwick, A. N. Norris,
\textit{Fundamental elastodynamic solutions for anisotropic
 media with ellipsoidal slowness surfaces},
Proc. Royal Soc. of London, \textbf{440(1910)}: (1993), 655--681.

\bibitem{Catheline} S. Catheline, J. L. Gennisson, G. Delon,  M. Fink, R. Sinkus, S. Abdouelkaram, J. Culioli,
\textit{Measuring of viscoelastic properties of homogeneous soft solid using transient elastography: an inverse problem approach},
Journal of Acoustical Society of America,
\textbf{116(6)}: (2004), 3734--3741.

\bibitem{Caputo} M. Caputo, 
\textit{Linear models of dissipation whose Q is
 almost frequency independent-II}, Geophysical  Journal  International, 
\textbf{13(5)}: (1967), 529–-539.

\bibitem{Carcoine2007Wave} J. M. Carcione,
\textit{Wave Field in the Real Media}, Elsevier Science, (second edition), 2007.


\bibitem{cerveny2008RayTheory} V. {\v{C}}erven{\`y},
\textit{Seismic Ray Theory},
Cambridge University Press, 2001.


\bibitem{Chandrasekhare} S. Chandrasekhar,
\textit{Ellipsoidal Figures of Equilibrium}, 
Yale University Press, 1969.


\bibitem{CourantHilbert} R. Courant, D. Hilbert,
\textit{Methods of Mathematical Physics}, 
\textbf{Vol 2}, Wiley-Interscience, 1989.


\bibitem{Dellinger} J. Dellinger,
\textit{Anisotropic Seismic Wave Propagation},
PhD Thesis, Stanford University, 1991.

\bibitem{GenCathChafFink} J. L. Gennisson, S. Catheline,
 S. Chaffai, M. Fink,
\textit{Transient elastography in anisotropic medium: 
Application to the measurement 
of slow and fast shear wave speeds in muscles},
J. Acous. Soc. Am., \textbf{114(1)}: (2003), 536--541.

\bibitem{Greenleaf} J.F. Greenleaf, M. Fatemi, M. Insana,
\textit{Selected methods for imaging elastic properties of biological tissues},
Annu. Rev. Biomed. Eng., \textbf{5}: (2003), 57--78.

\bibitem{Helbig} K. Helbig, 
\textit{Foundations of Anisotropy for Exploration Seismics},
Pergamon, New York, 1994.

\bibitem{HelbigThomsen} K. Helbig, L. Thomsen,
\textit{75-plus years of anisotropy in exploration and reservoir
seismics: A historical review of concepts and methods},
Geophysics, \textbf{70(6)}: (2006).

\bibitem{kellogg1953foun} O. D. Kellogg,
\textit{Foundations of Potential Theory},
Frederick Unger Publishing Company, New York, 1929.
 

\bibitem{lekhnitskii1964theory} S. G. Lekhnitskii,
\textit{Theory of Elasticity of an Anisotropic Body},
Mir Publishers, Moscow, 1981.


\bibitem{NamaniBaylay09} R. Namani, P. V. Bayly, 
\textit{Shear wave propagation in anisotropic soft tissues and gels},
EMBC 2009. Annual International Conference of IEEE, Minneapolis, USA, (2009).


\bibitem{Nedelec} J. C. N\'ed\'elec,
\textit{Acoustic and Electromagnetic Equations}, Applied Mathematical Sciences, \textbf{vol. 144}, Springer Verlag, 2001.

\bibitem{OidaKang} T. Oida, Y.Kang, T. Azuma, J. Okamoto,
 A. Amano, L. Axel, O. Takizawa, S. Tsutsumi, T. Matsuda,
\textit{The measurement of anisotropic elasticity in skeletal
 muscle using MR Elastography},
Proc. Intl. Soc. Mag. Reson. Med., 
\textbf{13}: (2005), 2020--2020.

\bibitem{Payton} R. G. Payton,
\textit{Elastic Wave Propagation in Transversely Isotropic Media},
Martinus Nijhoff Publishers, 1983.

\bibitem{Sarvazyan}  A. P. Sarvazyan, O. V. Rudenko, S. C. Swanson, J.B. Fowlkers, S. V. Emelianovs,
\textit{Shear wave elasticity imaging: a new ultrasonic technology of medical diagnostics},
Ultrasound in Med. \& Biol., \textbf{24(9)}: (1998), 1419--1435.

\bibitem{Sinkus00} R. Sinkus, J. Lorenzen, D. Schrader,
 M. Lorenzen, M. Dargatz, D. Holz,
\textit{High resolution tensor MR Elastography for breast tumor detection},
Phys. Med. Biol., \textit{45}: (2000).

\bibitem{Sinkus05} R. Sinkus, M. Tanter, S. Catheline, J. Lorenzen, C. Kuhl, 
E. Sondermann, M. Fink,
\textit{Imaging anisotropic and viscous properties of breast tissue by magnetic resonance-elastography}, Magnetic Resonance in Medicine, 
\textbf{53(2)}:(2005), 372--387.

\bibitem {SzaboWu00} T. L. Szabo, J. Wu, 
\textit{A model for longitudinal and shear wave propagation
in viscoelastic media}, 
Journal of Acoustical Society of America,
\textbf{107(5)}: (2000), 2437--2446.

\bibitem {SzaboCausal94} T. L. Szabo, 
\textit{Time domain wave equations for lossy media obeying
a frequency power law},
Journal of Acoustical Society of America,
\textbf{96(1)} (1994), 491--500.


\bibitem{titchmarsh} E. C. Titchmarsh,
\textit{Introduction to the Theory of Fourier Integrals},
Clarendon Press Oxford, (second edition) 1948.


\bibitem {VavAsymtGreen07} V. Vavry\u{c}uk, 
\textit{Asymptotic Green's function in homogeneous anisotropic
 viscoelastic media},
Proc. Royal Soc. A,
\textbf{463}: (2007), 2689--2707.

\bibitem{vavrycuk2001exact} V. Vavry\u{c}uk,
\textit{Exact elastodynamic Green functions for simple types 
of anisotropy derived from higher-order ray theory},
Studia Geophysica \& Geodaetica,
\textbf{45(1)}: (2001), {67--84}.

\bibitem{vavrycukWeak} V. Vavry\u{c}uk,
\textit{Elastodynamic and elastostatic Green tensors for homogeneous
weak transversely isotropic media},
Geophysics J. Int.,
\textbf{130}: (1997), {786--800}.

\bibitem{Weaver} J. Weaver, M. Doyley, E. Van Houten,
 M. Hood, X. C. Qin, F. Kennedy, S. Poplack, K. Paulsen,
\textit{Evidence of the anisotropic nature of the mechanical
 properties of breast tissue.}
 Med. Phys., \textbf{29}: (2002), 1291--1291

\bibitem{YoonKatz} H. S. Yoon, J. L. Katz,
\textit{Ultrasonic wave propagation in human cortical bone-I. Theoretical 
considerations for hexagonal symmetry},
J. Biomech., \textbf{9(6)}: (1976), 407--412.

%\bibitem{YapingZhu} Y. Zhu,
%\textit{Seismic Wave Propagation in
%Attenuative Anisotropic Media},
%PhD Thesis. Colorado School of Mines, 2006.

\end{thebibliography}

\end{document}